   \documentclass[11pt,letterpaper]{amsart}
\usepackage{graphicx}
\usepackage{fancyhdr}
\usepackage{amsmath, amssymb,amsthm}
\usepackage[all]{xy}
\usepackage{float}
\usepackage{tikz}
\usepackage{pdflscape}
\usepackage{longtable}
\usepackage{rotating}
\usepackage{hyperref}
\usepackage{subfigure}
\usepackage{tensor}
\usepackage{caption}
\usepackage{stmaryrd}
\usepackage{enumerate}

\usepackage{color}
\theoremstyle{plain}
\newtheorem{thm}{Theorem}

\newtheorem{lem}[thm]{Lemma}

\theoremstyle{definition}

\theoremstyle{remark}
\newtheorem{rem}[thm]{Remark}
\numberwithin{equation}{section}

\newcommand{\lgw}{\longrightarrow}

\newcommand{\si}{\sigma}

\newcommand{\ovl}{\overline}

\renewcommand{\Im}{\operatorname{Im}}

\newcommand{\la}{\lambda}

\renewcommand{\geq}{\geqslant}

\renewcommand{\k}{\Bbbk}

\newcommand{\R}{\mathbb{R}}

\newcommand{\K}{\mathbb{K}}

\newcommand{\N}{\mathbb{N}}

\newcommand{\C}{\mathbb{C}}
\newcommand{\Q}{\mathbb{Q}}
\newcommand{\LL}{\mathbb{L}}

\renewcommand{\lg}{\langle}
\newcommand{\rg}{\rangle}
\newcommand{\lb}{\llbracket}
\newcommand{\rb}{\rrbracket}

\renewcommand{\a}{\alpha}

\renewcommand{\b}{\beta}
\newcommand{\g}{\gamma}
\renewcommand{\phi}{\varphi}
\renewcommand{\d}{\delta}

\newcommand{\e}{\varepsilon}

\DeclareMathOperator{\im}{Im}

\begin{document}
\title[Not algebraically equivalent real surfaces]{Real algebraic surfaces biholomorphically equivalent but not algebraically equivalent}
\author{Guillaume Rond}
\email{guillaume.rond@univ-amu.fr}
\address{Aix Marseille Universit\'e, CNRS, I2M, Marseille, France}

\begin{abstract}
 We answer in the negative the long-standing open question of whether biholomorphic equivalence implies algebraic equivalence for germs of real algebraic manifolds in $\C^n$. More precisely we give an example of two germs of real algebraic surfaces in $\C^2$ that are biholomorphic, but not by means of  an algebraic biholomorphism.
\end{abstract}

\thanks{This research was funded, in whole or in part, by l'Agence Nationale de la Recherche (ANR), project ANR-22-CE40-0014. For the purpose of open access, the author has applied a CC-BY public copyright licence to any Author Accepted Manuscript (AAM) version arising from this submission.}


\maketitle

\section{Introduction}
Given two germs of smooth real analytic manifolds $(M,0)$ and $(M',0)$ in $\C^n$, we say that $(M,0)$ and $(M',0)$ are \emph{(biholomorphically) equivalent} if there is a germ of biholomorphism $\Phi:(\C^n,0)\lgw (\C^n,0)$ such that $\Phi(M)=M'$. The classification of real analytic manifolds up to local biholomorphisms is an old and important problem that goes back to H. Poincar\'e \cite{Po}, when he showed that  real analytic hypersurfaces of $\C^2$ have local invariants.  E. Cartan, for germs of real analytic smooth hypersurfaces in $\C^2$ \cite{Ca}, and S. S. Chern and J. K. Moser, for germs of Levi non-degenerate real analytic smooth hypersurfaces in $\C^n$ for $n\geq 2$ \cite{CM}, gave a complete description of this classification. More precisely, S. S. Chern and J. K. Moser first gave a complete classification up to formal biholomorphisms. Then they prove that the unique formal biholomorphism sending such a Levi non-degenerate hypersurface to its normal form  is convergent.\\
Therefore a natural question was to investigate if the formal  equivalence implies the convergent  equivalence. This question has been widely studied and the reader can consult \cite{Mir} for a general account of this problem. The first negative answer to this question has been given in \cite{MW}: the authors considered a particular example of two  real algebraic smooth surface germs $(M,0)$ and $(M',0)$ that are formally equivalent but not biholomorphically equivalent. This surface $M$ has the following particular property: its tangent space at  any point near the origin is totally real, but its tangent space at the origin is a complex line. A surface having this property is called a Bishop surface. \\
Recently this question has also been answered in the negative in \cite{KS} for CR manifolds, that is, for manifolds for which the largest $\C$-vector subspace of its tangent space has constant dimension.\\
\\
Here we consider the case of real algebraic manifolds in $\C^n$. One can define the notion of  \emph{algebraic  equivalence}:  two germs of smooth real algebraic manifolds $(M,0)$ and $(M',0)$ in $\C^n$ are algebraically equivalent if there is a germ of biholomorphism $\Phi:(\C^n,0)\lgw (\C^n,0)$, defined by algebraic power series, such that $\Phi(M)=M'$.  A formal power series $f(x_1,\ldots, x_n)\in\C\lb x_1,\ldots, x_n\rb$ is called \emph{algebraic} if it satisfies a non trivial relation
$$a_0(x)f(x)^d+a_1(x)f(x)^{d-1}\cdots+a_d(x)=0$$
where the $a_i(x)$ are polynomials. The question of whether biholomorphic equivalence implies  algebraic equivalence of germs of real algebraic manifolds in whole generality has first been asked in \cite[7. Question (b)]{BER} (see also \cite[Question B]{Mir}). But this question had  already been investigated before: H. Poincar\'e  studied algebraicity properties of local biholomorphisms between real algebraic hypersurfaces \cite{Po}, and he proved that  local biholomorphisms between  pieces of 3-spheres in $\C^2$ are necessarily  rational mappings (this has been extended in higher dimension by Tanaka \cite{Ta}).  Then an important step has been made by S. M. Webster who proved that
biholomorphisms between Levi non-degenerate real algebraic hypersurfaces  are necessarily 
algebraic \cite{Web1}. Now  the answer is known to be positive in several cases: for example the case of real algebraic hypersurfaces \cite{Web1}, \cite{BMR}, or  the case of real algebraic generic manifolds of finite type holomorphically non-degenerate \cite{BER0}. It is known that, for CR manifolds, the biholomorphic equivalence implies  the algebraic equivalence on a Zariski dense subset of points \cite{BRZ, LM}. See also \cite{Suk, Huang, BR, SS, CMS,Mer, Mi2, KLS,Z} for other related results, extensions and references, and \cite{Mir} for a survey about this question.\\
\\
In this paper we give an example of two algebraic Bishop surfaces in $\C^2$ that are biholomorphic but not algebraically biholomorphic, showing that the answer to this question is negative in general. Such surfaces have first been studied by E. Bishop in \cite{Bi} where he proved that they are locally biholomorphic to a surface defined (locally at 0) by:
$$w=z\ovl z +\g (z^2+\ovl z^2) +O(|z|^3) \text{ and } \im(w)=0.$$
The constant $\g$ is a biholomorphic invariant of the germ $(M,p)$, and is called the \emph{Bishop invariant} of $(M,p)$. J. K. Moser and S. M. Webster proved that, for $\g\notin\{0,\frac12,\infty\}$, a Bishop surface admits a (formal) normal form as follows:
$$w=z\ovl z +(\g+\e w^s) (z^2+\ovl z^2)  \text{ and } \im(w)=0$$
where $\e\in\{-1,0,1\}$ and $s\in\N^*\cup\{\infty\}$. 
 They also proved that this normal form can be obtained by a convergent biholomorphism when $\g\in\ ]0,\frac12[$ (see \cite[Theorem 1]{MW}). 
Here we prove the following:

\begin{thm}\label{main:thm}
Let $M_\g$ be the Bishop surface defined by
$$M_\g:=\{(z,w)\in\C^2\mid w=|z|^2+\g(\ovl z^2+z^3\ovl z)\}.$$
Assume that $\g\in\ ]0,\frac12[$ is a transcendental number. Then $M$ is biholomorphically equivalent to its normal form, but not by means of an algebraically biholomorphic map.
\end{thm}

 We  thank Nordine Mir who pointed out a mistake in a
previous version, and his multiple helpful remarks and valuable suggestions on this work. We also thank Andrew Elvey Price and Kilian Raschel for the fruitful discussions we
had on this problem.

\section{Algebraic power series}
In this part we review some results about algebraic power series that we will use in the proof of Theorem \ref{main:thm}. The reader can consul \cite{Ro0} and the references therein for more details about the ring of algebraic power series.

Let $\K$ be a field and $x$ denote the vector of indeterminates $(x_1,\ldots, x_n)$. A formal power series $f(x_1,\ldots, x_n)\in\K\lb x_1,\ldots, x_n\rb$ is called \emph{algebraic} if it satisfies a non trivial relation
$$a_0(x)f(x)^d+a_1(x)f(x)^{d-1}\cdots+a_d(x)=0$$
where the $a_i(x)$ are polynomials. The set of algebraic power series forms a subring of $\K\lb x\rb$ denoted by $\K\lg x\rg$. This is a Noetherian local ring, and, when $\K=\C$ every algebraic power series is a convergent power series. This ring satisfies the Implicit Function Theorem in the following sense: let $x=(x_1,\ldots, x_n)$ and $y=(y_1,\ldots, y_m)$ be two tuples of variables. Let $f_1(x,y)$, \ldots, $f_m(x,y)\in\K\lg x,y\rg$ such that
$$f_1(0,0)=\cdots=f_m(0,0),\text{ and the matrix } \frac{\partial f}{\partial y}(0,0) \text{ is invertible.}$$
Then there is a unique vector of algebraic power series $y(x)=(y_1(x),\ldots, y_m(x))\in\K\lg x\rg^m$ such that
$$f_1(x,y(x))=\cdots=f_m(x,y(x)) \text{ and } y_1(0)=\cdots=y_m(0)=0.$$
In fact  $\K\lg x\rg$ is the smallest subring of $\K\lb x\rb$ containing $\K[x]$ and satisfying the Implicit Function Theorem. \\
We mention that the composition of algebraic power series is again an algebraic power series. Moreover if a power series $f(x)$ with complex coefficients is algebraic, then its conjugate $\ovl f(x)$ is also algebraic.  

The following remarks will be useful in the sequel:
\begin{rem}\label{rem:diff alg}
Let $\Phi:(\K^n,0)\lgw (\K^n,0)$ be a bi-analytic automorphism induced by algebraic power series. Then the inverse $\Phi^{-1}$ is also induced by algebraic power series. This comes from the Inverse Function Theorem that is equivalent to the Implicit Function Theorem.
\end{rem}

\begin{rem}\label{flatness}
Since the ring $\K\lg x\rg$ is a Noetherian local ring, the ring extension $\K\lg x\rg\lgw \K\lb x\rb$ is faithfully flat (see \cite[Theorem 8.14]{Mat} for example). This implies that any system of linear equations
$$a_{j,1}(x)y_1+\cdots+a_{j,p}y_p=b_j(x), \text{ for }j=1, \ldots, q,$$
where the $a_{j,k}(x)$ and $b_j(x)\in\K\lg x\rg$, that has a formal power series solution 
$$(f_1(x),\ldots, f_p(x))\in\K\lb x\rb^p,$$ has also a solution in $\K\lg x\rg^p$ (see \cite[Example 1.4]{RoSurvey} for example).
\end{rem}

\begin{rem}\label{extension_field_power_series}
Let $\K\lgw \LL$ be an algebraic field extension. Then $\K\lb x\rb\cap\LL\lg x\rg=\K\lg x\rg$. Indeed, if $f(x)\in\K\lb x\rb\cap\LL\lg x\rg$, then $f(x)$ is algebraic over $\LL[x]$. But $\K[x]\lgw \LL[x]$ is an algebraic extension of rings. Therefore, $f(x)$ is algebraic over $\K[x]$.
\end{rem}

\begin{lem}\label{finite-extension-series}
Let $\K$ be field and $c$ be finite over $\K$ of degree $d>1$. Assume that $c$ is separable over $\K$. Let $x=(x_1,\ldots, x_n)$ be a tuple of indeterminates, and let $\phi(x)\in\K[c]\langle x\rangle$. Let us write
$$\phi(x)=\phi_0(x)+c\phi_1(x)+\cdots+c^{d-1}\phi_{d-1}(x)$$
where the $\phi_j$ are formal power series with coefficients in $\K$. Then, for every $j$, $\phi_j(x)\in\K\lg x\rg$.
\end{lem}

\begin{proof}
The power series $\phi(x)$ is algebraic over $\K[c][x]$. But $\K[c][x]$ is a finite extension of $\K[x]$, thus $\phi(x)$ is algebraic over $\K[x]$. Let $P(x,y)\in\K[x,y]$ be nonzero with $P(x,\phi(x))=0$.

Let us denote by $c_1:=c$, $c_2$, \ldots, $c_d$ the conjugates of $c$ over $\K$ in a given algebraic closure $\ovl \K$ of $\K$. Let us set, for $j=1$, \ldots, $d$
$$\phi(c_j,x):=\phi_0(x)+c_j\phi_1(x)+\cdots+c_j^{d-1}\phi_{d-1}(x).$$
The $\phi(c_j,x)$ are algebraic power series since $P(x,\phi(c_j,x))=0$ for every $j$. Indeed, if $\si$ is a $\K$-automorphism of $\ovl\K$ sending $c=c_1$ on $c_j$,  we have
$$0=P(x,\phi(c,x))=\si(P(x,\phi(c,x))=P(x,\phi(\si(c),x))=P(x,\phi(c_j,x))$$
since the coefficients of $P$ belong to $\K$.
We have 
$$\begin{pmatrix} 1 & c_1 & c_1^2 & \cdots & c_1^{d-1}\\
1 & c_2 & c_2^2 & \cdots & c_2^{d-1}\\
\vdots & \vdots & \vdots & \cdots & \vdots\\
1 & c_d & c_d^2 & \cdots & c_d^{d-1}
\end{pmatrix}
\begin{pmatrix}\phi_0(x) \\ \phi_1(x) \\ \vdots\\ \phi_{d-1}(x)
\end{pmatrix}=
\begin{pmatrix}\phi(c_1,x) \\ \phi(c_2,x)\\ \vdots \\ \phi(c_d,x)
\end{pmatrix}$$
Let us denote by $M$ the left-hand Vandermonde matrix. This matrix is invertible and the entries of its inverse belong to $\ovl \K$. Therefore, by Remark \ref{extension_field_power_series}, $\phi_0(x)$, \ldots, $\phi_{d-1}(x)$ belong to $\K\lg x\rg$.
\end{proof}

\begin{lem}\label{Gilmer}
Let $x=(x_1,\ldots, x_n)$, $\K$ be a characteristic zero field, $\k$ be a subfield of $\K$ and  $f=\sum_{\a\in\N^n}f_\a x^\a \in\K\langle x\rangle$. Then there is a field $\LL\subset \K$, finitely generated over $\k$, such that
$$\forall \a\in \N^n,\ f_\a\in\LL.$$
\end{lem}
\begin{proof}
This result is well known but we did not find an explicit reference. This result is an easy corollary of \cite{G}:
let $\K_0$ be the  field extension of $\k$ generated by the coefficients of a nonzero polynomial $P(x,y)$ such that $P(x,f(x))=0$. So $\K_0$ is a  finitely generated field extension of $\k$. Moreover $f$ is algebraic over $\K_0[x]$, so by \cite[Theorem 1.1]{G}, all the $f_\a$ are algebraic over $\K_0$. Then by \cite[Theorem 2.4]{G} or \cite[Theorem 1.1]{CK}, the $f_\a$ belong to a finite field extension $\LL$ of $\K_0$. In particular $\K_0\lgw \LL$ is a finitely generated field extension. So $\k\lgw \LL$ is finitely generated.
\end{proof}

\begin{lem}[Eisenstein's Theorem]\label{key-lemma} (see \cite[Lemma 11]{RoANT} or \cite{Ro24} for example)
Let $x=(x_1,\ldots, x_n)$ and $t=(t_1,\ldots,t_r)$ be two tuples of indeterminates and $\K$  be a field extension of $\ovl\Q$ generated by algebraically independent elements $\g_1$, \ldots, $\g_r$.
Let $f\in\K\langle x\rangle^m$ be a tuple of algebraic power series. Then there exist a polynomial $S(t)\in\ovl\Q[t]$ and, for every $\a\in\N^n$, polynomials $R_{j,\a}(t)\in\ovl\Q[t]$ such that
$$\forall j=1,\ldots, m,\ f_j(x)=F_j(\g_1,\ldots,\g_r,x) \text{ with } F_j(t,x)=\sum_{\a\in\N^n}\frac{R_{j,\a}(t)}{S(t)^{|\a|}}x^\a$$
where $|\a|:=\a_1+\cdots+\a_n$.
\end{lem}

\section{Proof of the main theorem}
Let $\g\in\ ]0,1/2[$. By \cite[Theorem 1]{MW}, there is a convergent biholomorphic map germ $\Phi:(\C^2,0)\lgw (\C^2,0)$ such that $\Phi(M_\g)$ is the germ of some algebraic set $N_{\g,\e,s}$ as defined in the introduction. In particular there is a convergent biholomorphic map $\Phi$ sending $M_\g$ into the real hyperplane $\Im(w)=0$. If we denote by $\phi_1(z,w)$ and $\phi_2(z,w)\in\C\{z,w\}$ the components of $\Phi$, this means that 
$$\phi_2(z,w)-\ovl\phi_2(\ovl z,\ovl w)=0$$
when $w=|z|^2+\g(\ovl z^2+z^3\ovl z)$. This is equivalent to
$$\phi_2(z,|z|^2+\g(\ovl z^2+z^3\ovl z))-\ovl\phi_2(\ovl z,\ovl{|z|^2+\g(\ovl z^2+z^3\ovl z)})=0.$$
This is also equivalent to the existence of two power series
 $\ell_1(z,\ovl z,w,\ovl w)$, $\ell_2(z,\ovl z,w,\ovl w)\in\C\{z,\ovl z,w,\ovl w\}$ such that
\begin{equation}\label{main_equation}
\begin{split}
\phi_2(z,w)-\ovl\phi_2(\ovl z,\ovl w)=\ell_1(z,\ovl z,w,\ovl w)&\left(w-|z|^2-\g(\ovl z^2+z^3\ovl z)\right)\\
&+\ell_2(z,\ovl z,w,\ovl w)\left(\ovl w-|z|^2-\g(z^2+\ovl z^3 z)\right)\end{split}
\end{equation}

When $\g>1/2$ is such that the solutions of the equation $\g T^2-T+\g=0$ are not roots of unity (one says that $\g$ is \emph{exceptional}), Moser and Webster proved that there is a formal biholomorphic map germ $\Phi$ satisfying Equation \eqref{main_equation}, but that this formal map is not convergent (see  \cite[Theorem 6.1]{MW}). 

Now let us fix $\g_1\in\ ]0,1/2[$ to be a transcendental number. The proof is done by contradiction : assume that there exists  a convergent biholomorphic map germ $\Phi:(\C^2,0)\lgw (\C^2,0)$ whose components, $\phi_1(z,w)$, $\phi_2(z,w)\in\C\lg z,w\rg$ are algebraic power series and satisfy \eqref{main_equation}. 
Since $\phi_2(z,w)$, $\ovl\phi_2(\ovl z,\ovl w)$, $w-|z|^2-\g_1(\ovl z^2+z^3\ovl z)$ and $\ovl w-|z|^2-\g( z^2+\ovl z^3z)$ are algebraic power series, by Remark \ref{flatness} we may assume that $\ell_k(z,\ovl z, w, \ovl w)\in\C\lg z,\ovl z,w,\ovl w\rg$ for $k=1$, $2$. 

By Lemma \ref{Gilmer} the coefficients of the algebraic series $\phi_2(z,w)$ and the $\ell_k(z,\ovl z, w, \ovl w)$, for $k=1$, $2$,  belong to a finitely generated field extension $\LL$ of $\K:=\ovl \Q(\g_1)$. Let $\xi_1$, \ldots, $\xi_p\in\C$ be generators of $\LL$: $\LL=\K(\xi_1,\ldots, \xi_p)$. Since $i\in\ovl\Q\subset \K$, we may replace the $\xi_k$ by their real and imaginary parts  in order to assume that $\LL$ is a finitely generated field extension of $\K$ generated by real numbers. Take a maximal algebraically independent subfamily of $\{\xi_1,\ldots, \xi_p\}$ and denote the elements of this family by $\g_2$, \ldots, $\g_r$. Thus $\LL$ is a finite field extension of $\LL_1:=\K(\g_2,\ldots,\g_r)$. By the primitive element theorem, this finite extension is generated by one algebraic element that we  write as $a+ib$ where $a$, $b\in\R$. So the coefficients of the algebraic series are in $\LL_1[a,b]$. But, it is well known that $\LL_1[a,b]=\LL_1[a+\la b]$ for all but finitely many $\la\in\LL_1$. So we may find $\la\in\Q$ such that $c:=a+\la b\in\R $ is finite over $\LL_1$ and the coefficients of the algebraic series belong to $\K(\g_2,\ldots,\g_r)[c]$.

If $d$ denotes the degree of $c$ over $\LL_1$, we can write
$$\phi_2(z,w)=\sum_{j=0}^{d-1}c^j\phi_{2,j}(z,w)\text{ and } \ell_k(z,\ovl z,w,\ovl w)=\sum_{j=0}^{d-1}c^j\ell_{k,j}(z,\ovl z,w,\ovl w)$$
where the $\phi_{2,j}$ and the $\ell_{k,j}$ are power series with coefficients in $\LL_1$. By Lemma \ref{finite-extension-series}, the $\phi_{2,j}$ and the $\ell_{k,j}$ are algebraic power series with coefficients in $\LL_1$. Moreover, because  \eqref{main_equation} is a linear equation, because $c\in\R$, and because the coefficients of $w-|z|^2-\g(\ovl z^2+z^3\ovl z)$ and $\ovl w-|z|^2-\g(z^2+\ovl z^3 z)$ are in $\LL_1$, 
$(\phi_2,\ell_1,\ell_2)$ is a solution of \eqref{main_equation} if and only if all the $(\phi_{2,j},\ell_{1,j},\ell_{2,j})$ are solutions of \eqref{main_equation}. Moreover, since $\Phi$ is a germ of biholomorphism, the vanishing order of $\phi_2(z,w)$ is 1. So the vanishing order of  $\phi_{2,j}$, for some $j$ (let us say $j_0$), is also 1, since $1$, $c$, \ldots, $c^{d-1}$ are $\LL_1$-linearly independent. Therefore we may replace $(\phi_2,\ell_1,\ell_2)$ by $(\phi_{2,j_0},\ell_{1,j_0},\ell_{2,j_0})$ and assume that the coefficients of $\phi_2$ and the $\ell_k$ are in $\LL_1$.

Now, by Lemma \ref{key-lemma}, we can write
\begin{equation}\label{eq1}\phi_2(z,w)=\sum_{\a\in\N^2}\frac{R_\a(\g_1,\ldots, \g_r)}{S(\g_1,\ldots,\g_r)^{|\a|}}z^{\a_1}w^{\a_2}\end{equation}
\begin{equation}\label{eq2}\text{ and }  \ell_k(z,\ovl z,w,\ovl w)=\sum_{(\a,\b)\in\N^4}\frac{R_{k,\a,\b}(\g_1,\ldots, \g_r)}{S(\g_1,\ldots,\g_r)^{|\a|+|\b|}}z^{\a_1}w^{\a_2}\ovl z^{\b_1}\ovl w^{\b_2} \end{equation}
for some polynomials $R_{\a}$, $S$, $R_{k,\a,\b}\in\ovl \Q[t_1,\ldots, t_r]$. 
 We set $t=(t_1,\ldots, t_r)$ and 
 $$U(t,z,w):=\sum_{\a\in\N^2}\frac{R_\a(t_1,\ldots, t_r)}{S(t_1,\ldots,t_r)^{|\a|}}z^{\a_1}w^{\a_2}$$
 and
 $$V_k(t,z,\ovl z,w,\ovl w)=\sum_{(\a,\b)\in\N^4}\frac{R_{k,\a,\b}(t_1,\ldots, t_r)}{S(t_1,\ldots,t_r)^{|\a|+|\b|}}z^{\a_1}w^{\a_2}\ovl z^{\b_1}\ovl w^{\b_2}.$$ 
 Let $P(t,z,w,u)\in\ovl \Q(t)[z,w,u]$ and $Q_k(t,z,\ovl z,w,\ovl w,v)\in\ovl\Q(t)[z,\ovl z,w,\ovl w,v]$  be irreducible polynomials such that
 \begin{equation}\label{EE}P(\g_1,\ldots,\g_r,z,w,\phi_2(z,w))=0,\ Q_k(\g_1,\ldots,\g_r,z,\ovl z,w,\ovl w , \ell_k(z,\ovl z,w,\ovl w))=0.\end{equation}
 Up to multiplication by an element of $\ovl\Q(t)$ we may assume that each of the three polynomials $P$, $Q_1$ and $Q_2$ has a monomial (in the indeterminates $z$, $\ovl z$, $w$, $\ovl w$, $u$ and $v$) whose coefficient in $\ovl\Q(t)$ equals 1.
 We have that
 \begin{equation}\label{EEE}P(t,z,w,U(t,z,w)),\ Q_k(t,z,\ovl z,w,\ovl w,V_k(t,z,\ovl z,w,\ovl w))\in\ovl\Q(t)\lb z,\ovl z,w,\ovl w\rb\end{equation}
 and, by \eqref{EE}, the coefficients of these three power series vanish when we set $t=\underline\g=(\g_1,\ldots,\g_r)$. The coefficients of these three power series are rational functions whose denominators are products of a finite number of polynomials in $t$: the polynomial $S(t)$ and the (finitely many) denominators of the coefficients of $P$, $Q_1$ and $Q_2$. The numerators of the three series in \eqref{EEE} are polynomials in $t_1$, \ldots, $t_r$ that vanish at $\underline\g$. Since the $\g_i$ are algebraically independent over $\ovl\Q$, these numerators are identically zero. So these series in \eqref{EEE} are identically zero.

 Let $\mathcal V$ be the common zero locus in $\R^r$ of $S(t)$ and of the denominators of the coefficients of $P$, $Q_1$ and $Q_2$. The set $\mathcal V$ is a proper algebraic subset of $\R^r$ and $\underline\g\in\R^r\setminus \mathcal V$. 
 And for $\underline\d\in\R^r\setminus \mathcal V$, the series $U(\underline\d,z,w)$, the $V(\underline\d, z,\ovl z,w,\ovl w)$  and the polynomials $P(\underline\d,z,w,u)$ and  $Q(\underline\d,z,\ovl z,w,\ovl w,v)$ are well defined. Moreover, because the series in \eqref{EEE} are zero, $U(\underline\d,z,w)$ and the $V_k(\underline\d, z,\ovl z,w,\ovl w)$ are algebraic power series since they are roots of $P(\underline\d,z,w,u)$ and of the $Q_k(\underline\d,z,\ovl z,w,\ovl w,v)$ (which are nonzero polynomials because $P$, $Q_1$ and $Q_2$ have a coefficient equal to 1). We can expand the following power series:
 \begin{equation*}\begin{split}&U(t,z,w)-\ovl U(t,\ovl z,\ovl w)-V_1(t, z,\ovl z,w,\ovl w)\left(w-|z|^2-t_1(\ovl z^2+z^3\ovl z)\right)\\
 -&V_2(t, z,\ovl z,w,\ovl w)\left(\ovl w-|z|^2-t_1(z^2+\ovl z^3z)\right)=\sum_{\b\in\N^4}F_\b(t_1,\ldots, t_r)z^{\b_1}\ovl z^{\b_2}w^{\b_3}\ovl w^{\b_4}\end{split}\end{equation*}
 where the $F_\b(t)$ are rational functions in $\ovl\Q(t_1,\ldots,t_r)$, whose denominators do not vanish on $\R^r\setminus\mathcal V$. Since $(\phi_2,\ell_1,\ell_2)$ is a solution of \eqref{main_equation}, we have that
 $$\forall \b\in\N^4,\ F_\b(\underline \g)=0.$$
 Since $\g_1$, \ldots, $\g_r$ are algebraically independent over $\ovl\Q$, we have that, for every $\b$, $F_\b(t)=0$.
 
 \begin{lem}
 There is a finite set $E\subset\ovl\Q$ such that
 $$\forall \d_1\in\R\setminus E,\ \exists \d_2, \ldots, \d_r \in\R \text{ such that } (\d_1,\ldots,\d_r)\in\R^r\setminus \mathcal V.$$
 \end{lem}
 
 \begin{proof}
 Let $\d_1\in\R$. Assume that  for all $\d_2$, \ldots, $\d_r\in\R$, $(\d_1,\ldots,\d_r)\in \mathcal V$. Let $A(t)$ be one polynomial among the following polynomials: $S(t)$ and the (finitely many) denominators of the coefficients of $P$ and  $Q$. By assumption $A(\d_1,t_2,\ldots, t_r)$ vanishes on $\R^{r-1}$. Thus $A(\d_1,t_2,\ldots, t_r)=0$, so $A(t)$ is divisible by $(t-\d_1)$. But because $\C[t_1,\ldots, t_r]$ is a unique factorization domain, there is only finitely many $\d_1\in\R$ such that $t-\d_1$ divides $A(t)$. We denote by $E$ the set of these elements. Moreover, since $A(t)\in\ovl\Q[t]$, $t-\d_1$ divides $A(t)$ only when $\d_1\in\ovl\Q$; indeed,  $t-\d_1$ divides $A(t)$ if and only if $\d_1$ is a root of every coefficient (in $\ovl\Q[t_1]$) of $A(t)$ seen as a polynomial in $t_2$, \ldots, $t_r$ $-$ but these coefficients are polynomials in $\ovl\Q[t_1]$, so their roots are in $\ovl\Q$.\\
Now if $t-\d_1$ does not divide $A(t)$, there exists $(\d_2,\ldots,\d_r)\in\R^{r-1}$ such that $A(\underline\d)\neq 0$, so $(\d_1,\ldots,\d_r)\in\R^r\setminus \mathcal V$.
 \end{proof}
 
 So finally, for  $\d_1\in\R\setminus E$, there is $\d_2$, \ldots, $\d_r \in\R$ such that
 $$\phi'_2(z,w):=U(\underline\d,z,w) \text{ and } \ell'_k(z,\ovl z,w,\ovl w):=V_k(\underline\d, z,\ovl z,w,\ovl w),\text{ for } k=1,2,$$
 are well defined algebraic power series satisfying
 \begin{equation}\label{main_equation2}
 \begin{split}\phi'_2(z,w)-\ovl\phi'_2(\ovl z,\ovl w)=\ell'_1(z,\ovl z,w,&\ovl w)\left(w-|z|^2-\d_1(\ovl z^2+z^3\ovl z)\right)\\
&+\ell_2(z,\ovl z,w,\ovl w)\left(\ovl w-|z|^2-\d_1(z^2+\ovl z^3 z)\right)\end{split}
\end{equation}
 and the linear term $L(z,w)$ of $\phi'_2(z,w)$ is nonzero. Let  $\phi'_1(z,w)$ be a homogeneous linear polynomial which is linearly independent with $L(z,w)$. Then the map
 $$\Phi':(\C^2,0)\lgw (\C^2,0),$$
 whose components are $\phi'_1$ and $\phi'_2$, is a local algebraic biholomorphic germ by Remark \Ref{rem:diff alg}. Since \eqref{main_equation2} is satisfied, we have that $\Phi'(M_{\d_1})\subset\{\Im(w)=0\}$. Moreover $\Phi'$ is a (convergent) biholomorphism since algebraic power series are convergent power series. But if we choose $\d_1>1/2$ to be exceptional (that is possible since $E$ is finite), such a $\Phi'$ cannot exists  by \cite[Theorem 6.1]{MW}. Therefore, by contradiction, $M_\g$ is biholomorphically equivalent to $N_{\g,\e,s}$ for some $(\e,s)\in\{1,0,-1\}\times(\N^*\cup\{\infty\})$, but these two germs are not equivalent by means of an algebraic biholomorphic map.

\end{document}